\documentclass[10pt,a4paper]{article}

\usepackage
[
        a4paper,
        left=3cm,
        right=3cm,
        top=3cm,
        bottom=4cm,
]
{geometry}

\usepackage{parskip}
\usepackage{amsthm}
\usepackage{amsmath}
\usepackage{amssymb}
\usepackage{xcolor}
\usepackage{natbib}
\usepackage{shuffle}   
\usepackage[utf8]{inputenc} 
\usepackage[T1]{fontenc}    
\usepackage{lmodern}
\usepackage{hyperref}       
\usepackage{url}            
\usepackage{booktabs}     
\usepackage{amsfonts}       
\usepackage{nicefrac}       
\usepackage{microtype}    
\usepackage{lipsum}
\usepackage{appendix}
\usepackage{forest}
\usepackage{enumitem}
\usepackage{tikz}
\usetikzlibrary{trees}
\usepackage{tikz-cd}
\usepackage{adjustbox}
\usepackage{algorithm}
\usepackage{float}
\usepackage[noend]{algpseudocode}
\usepackage{authblk}
\usepackage{cleveref}
\usepackage{dirtytalk}
\usepackage{cancel}

\theoremstyle{definition}
\newtheorem{theorem}{Theorem}[section]
\newtheorem{definition}[theorem]{Definition}

\newtheorem{lemma}[theorem]{Lemma}
\newtheorem{remark}[theorem]{Remark}
\newtheorem{example}[theorem]{Example}

\newcommand\DEF[1]{\emph{#1}}

\newlist{steps}{enumerate}{1}
\setlist[steps, 1]{label = Step \arabic*:}

\newcommand{\D}{\, \mathrm{d}}

\newcommand\area{\operatorname{area}}
\newcommand\vol{\operatorname{vol}}
\newcommand\areaTree{\mathsf{area}}

\newcommand\lhs{\mathbin{\prec}}
\newcommand\lhsf{\mathord{\prec}}

\newcommand{\Anonempty}{\mathcal{A}^{>0}}

\makeatletter
\renewcommand\AB@affilsepx{ - \protect\Affilfont}
\makeatother

\title{A structure theorem for streamed information}

\author[1,2]{\small Cristopher Salvi \footnote{Email: \url{c.salvi@imperial.ac.uk}}}
\author[3]{\small Joscha Diehl}
\author[2,4]{\small Terry Lyons}
\author[5]{\small Rosa Preiss}
\author[6]{\small Jeremy Reizenstein}

\affil[1]{\footnotesize Imperial College London}
\affil[2]{\footnotesize The Alan Turing Institute}
\affil[3]{\footnotesize University of Greifswald}
\affil[4]{\footnotesize University of Oxford}
\affil[5]{\footnotesize University of Potsdam}
\affil[6]{\footnotesize Meta AI}

\begin{document}

\maketitle

\begin{abstract}
    We identify the free half shuffle algebra of \citet{schutzenberger1958propriete} with an algebra of real-valued functionals on paths, where the half shuffle emulates integration of a functional against another. We then provide two, to our knowledge, new identities in arity 3 involving its commutator (area), and show that these are sufficient to recover the Zinbiel and Tortkara identities introduced by \citet{dzhumadil2007zinbiel}. We then use these identities to provide a simple proof of the main result of \citet{diehl2020areas}, namely that any element of the free half shuffle algebra can be expressed as a polynomial over iterated areas.
    
    Moreover, we consider minimal sets of Hall iterated integrals defined through the recursive application of the half shuffle product to 
    Hall trees. Leveraging the duality between this set of Hall integrals and classical Hall bases of the free Lie algebra, we prove using combinatorial arguments that any element of the free half shuffle algebra can be written uniquely as a polynomial over Hall integrals. We interpret this result as a structure theorem for streamed information, loosely analogous to the unique prime factorisation of integers, allowing to split any real valued function on streamed data into two parts: a first that extracts and packages the streamed information into recursively defined atomic objects (Hall integrals), and a second that evaluates a polynomial function in these objects without further reference to the original stream. The question of whether a similar result holds if Hall integrals are replaced by Hall areas is left as an open  conjecture.
    
    Finally, we construct a canonical, but to our knowledge, new decomposition of the free half shuffle algebra as shuffle power series in the greatest letter of the original alphabet with coefficients in a sub-algebra freely generated by a new alphabet with an infinite number of letters. We use this construction to provide a second proof of our structure theorem.
\end{abstract}




\section{Introduction}

It is not too much to accept that, at least on some fine enough time scales, most instance of streamed information (text, sound, video, time series...) can be represented,  as a path $\gamma : [0,1] \to V$ with values in some finite dimensional vector space $V \simeq \mathbb{R}^d$. It was first shown by \citet{chen1957integration}, and then explored in greater detail and generality in the context of \emph{rough path theory} in \citep{hambly2010uniqueness, boedihardjo2016signature}, that any path may be faithfully represented, up to reparameterisation, by the collection of its iterated integrals known as the \emph{signature}. This non-commutative exponential maps a path to a grouplike element on the tensor algebra $(\mathcal{A},\otimes)$, where $\mathcal{A}$ is the vector space spanned by words in $d$ letters, including the empty word $e$, and $\otimes$ is the tensor product. For an arbitrary interval $[a,b] \subset [0,1]$, the signature $\mathcal{S}(\gamma)_{a,b} := X_b$ where $X$ is the unique solution to the control system $dX_t = X_t \otimes d\gamma_t$ started at $X_a = e$. Furthermore, the range of the signature describes the set of characters $G \subset \mathcal{A}$.

The half shuffle product $\prec$ was firstly introduced in \citep{schutzenberger1958propriete}, where it also showed that $\mathcal{A}$ is the free algebra over $A$ with respect to $\prec$. We will later refer to this algebra as the \emph{free half shuffle algebra of Sch\"utzenberger}. In the same article, the shuffle product $\shuffle$ was subsequently defined as $f \shuffle g = f \prec g + g \prec f + \langle f, e\rangle \langle g, e\rangle e$, so to emulate integration by parts.

It is well known that the shuffle algebra $(\mathcal{A},\shuffle)$ is the algebraic dual of the tensor algebra $(\mathcal{A},\otimes)$ \citep{reutenauer1993london}; it is automatic from this perspective to see that the restriction of linear functionals on $\mathcal{A}$ to the range of the signature $G$ form a unital algebra of real-valued functions that separates points \citep{lyons2004differential}. A straightforward application of the Stone-Weierstrass theorem yields that for any compact set of reparameterisation-reduced paths, linear functionals acting on their signatures are dense in the space of continuous, real-valued functions on this compact set under a suitable choice of topology \citep{cass2022topologies}. 

Because $G$ is the set of characters, the main result in \citet{ree1958lie} implies that the restriction of the shuffle product of two of elements of the shuffle algebra to $G$ is the pointwise product of the two restrictions $\langle f \shuffle g, \mathcal{S}(\gamma)_{a,b}\rangle = \langle f, \mathcal{S}(\gamma)_{a,b}\rangle \langle g, \mathcal{S}(\gamma)_{a,b}\rangle$, the so-called \emph{shuffle identity}. This interplay between algebraic and analytic operations can be extended to the half shuffle product, emulating integration of a path functional against another $\langle f \lhs g, \mathcal{S}(\gamma)_{a,b}\rangle  = \int_a^b \langle g, \mathcal{S}(\gamma)_{a,s}\rangle d \langle f, \mathcal{S}(\gamma)_{a,s}\rangle$, and to its commutator representing the area enclosed by the two dimensional curve $ t \mapsto (\langle f, \mathcal{S}(\gamma)_{a,t}\rangle, \langle g, \mathcal{S}(\gamma)_{a,t}\rangle)$ and the chord connecting the two end points

$$\langle \area(f,g), \mathcal{S}(\gamma)_{a,b}\rangle = \int_a^b \langle g, \mathcal{S}(\gamma)_{a,s}\rangle d\langle f, \mathcal{S}(\gamma)_{a,s}\rangle - \int_a^b \langle f, \mathcal{S}(\gamma)_{a,s}\rangle d\langle g, \mathcal{S}(\gamma)_{a,s}\rangle.$$

Thus, collectively iterated integrals provide an accurate description of the path and linear combinations of them can be determined easily by regression, making the coefficient of the signature an ideal feature set for machine learning applications on streamed data \citep{fermanian2023new}; signature methods have been applied in a variety of contexts including deep learning for time series
\cite{kidger2019deep, morrill2021neural, cirone2023neural}, kernel methods \cite{salvi2021signature, lemercier2021distribution, lemercier2021siggpde} quantitative finance \cite{arribas2020sigsdes, salvi2021higher, horvath2023optimal} and cybersecurity \cite{cochrane2021sk}. 

However, these integrals contain some redundancies, in the sense that some higher ones can be expressed using polynomial relations in lower ones. This represents a major scalability issue, particularly because the number of distinct and linearly independent iterated integrals grows exponentially with the degree of iteration in the integral. This raises a simple set of questions which we will answer positively in this paper: 

\emph{Can we identify minimal sets of integrals so that each integral is an integral of two other integrals in the same class and so that every other integral can be expressed as a polynomial in them?} 

The minimal sets of integrals we identify in this paper are defined hierarchically using sets of binary planar rooted trees called \emph{Hall sets} \citep{reutenauer1993london, bourbaki2008lie}, and can be computed recursively in a localised way (to compute one, one must compute its ancestors but not others) which adds further value to the results. These minimal sets of integrals fully describes the information in the stream while the polynomials capture the nonlinearity in any function of interest. It is for this reason we call it a \emph{structure theorem}, loosely analogous to the unique factorisation of integers as products of primes. In this way we see that identifying a basis for the space of smooth functions acting on pathspace splits the evaluation process into two parts: a) a first that engages with the underlying stream of information\footnote{This information extraction is done in practice via some physical integration process that responds to the underlying signal. Physical integration processes are intrinsically nasty as mathematical operators (controlled differential equations in general, and in particular the integration process here, are not closable in the uniform topology on $\gamma$ - see \citep{lyons1998differential})}, systematically extracts and packages the relevant information into atomic objects whilst removing what's irrelevant, b) a second that evaluates a unique polynomial function in these expensive but informative precomputed basis elements in order to deliver the desired function evaluation without further reference to the original stream $\gamma$.

Having established that polynomials in \emph{Hall integrals} freely generate the half shuffle algebra $(\mathcal{A},\prec)$, it is natural to ask whether a similar structure theorem holds when the half shuffle $\prec$ is replaced by its commutator $\area$. This question has been, and still remain, a source of conjecture, well supported by calculation, for the last decade. Nonetheless, the search for an answer to this conjecture led us to consider an argument related to the well-known Lazard's elimination \citep{reutenauer1993london} to construct a canonical, but to our knowledge, new decomposition of the algebra $\mathcal{A}$ as shuffle power series in the greatest letter of the original alphabet with coefficients in a sub-algebra freely generated by a new alphabet with an infinite number of letters. This construction, that we refer to as elimination trick, will enable us to provide a second proof of our structure theorem relying on an induction argument.

We briefly outline the structure of the paper. Section \ref{sec:background} provides a brief background on the algebraic setup needed for the rest of the paper. In Section \ref{sec:identities} we introduce the free half shuffle algebra of Sch\"utzenberger, we make precise the interplay between the algebraic operations $\lhs,\shuffle,\area$ and the corresponding analytic operations on paths, and we provide two new identities in arity 3 involving the $\area$ product. In Section \ref{sec:iterated_areas} we make use of these new identities to provide a simpler proof of the main result in \citep{diehl2020areas}, stating that polynomials in iterated areas generate the algebra $\mathcal{A}$. In Section \ref{sec:structure_theorem} we present our structure theorem for streamed information, providing a simple proof of the main result in \citep{sussmann1986product} reported without proof also in \citep{kawski1999chronological,gehrig2008hopf} stating that polynomials in Hall integrals freely generate the algebra $\mathcal{A}$. Finally, using the elimination trick we provide a second proof of our structure theorem.

\section{Background}\label{sec:background}
First, we remind the reader in a very terse form of the general collection of objects about which we write. Much more can be found by looking in \citep{bourbaki2008lie} or (and we will follow this for the results we need) \citep{reutenauer1993london}.  We hope the paper is self contained, and cites what is needed, but for the rest of this introduction, we will be very brief and assume the reader has familiarity with the general algebraic framework. 

The starting point will be a finite alphabet $A$ of $d$ letters. 

\begin{definition}
    A \DEF{word} on the alphabet $A$ is a finite sequence of letters from $A$, including the empty sequence, called the \DEF{empty word} and denoted by $e$. We denote by $W_A$ the set of all words, including the empty word. $W_A$ with the concatenation product is a monoid, that is free over $A$. The length $|w|$ of a word $w \in W_A$ is the number of letters in $w$. Finally, we denote by $\mathcal{A}$ the vector space spanned by all words in $W_A$.
\end{definition}

\begin{remark}\label{rq:decomposition_A}
    The vector space $\mathcal{A}$ admits the unique direct sum decomposition 
    \begin{equation}\label{eqn:z_decomposition}
        \mathcal{A}= \Anonempty \oplus \langle e \rangle,
    \end{equation}
    where $\langle e \rangle$ is the vector space spanned by the empty word and $\Anonempty$ is its annihilator, i.e.
    $$\Anonempty := \{f \in \mathcal{A} : \langle f, e\rangle = 0\}.$$
    Note that $\Anonempty$ is the vector space spanned by all non-empty words.
    It follows that any $f \in \mathcal{A}$ admits the unique decomposition
    $$f = (f - \langle f, e\rangle e)  + \langle f, e\rangle e,$$
    where $(f - \langle f, e\rangle e) \in \Anonempty$ and $\langle f, e\rangle e \in \langle e \rangle$.
\end{remark}

$\mathcal{A}$ is graded by word length. The words of length greater than $n \in \mathbb{N}$ span an ideal, and the quotient of $\mathcal{A}$ by this ideal is often referred to as the \DEF{truncated tensor algebra} $\mathcal{A}^{(n)}$. 

\begin{definition}
Denote by $(\mathcal{A}, \otimes)$ the \DEF{tensor algebra over $A$}, that is the free associative $\mathbb{R}$-algebra over $A$ with the tensor product $\otimes$. 
\end{definition}

\begin{remark}\label{rq:pairing}
An infinite linear combination of words in $W_A$ is usually referred to as a \DEF{series}. There is a natural duality between $\mathcal{A}$ and the associative algebra of all series $\mathcal{A}^\infty$ given by the pairing $ (\cdot,\cdot) : \mathcal{A} \times \mathcal{A}^\infty \to \mathbb{R}$ defined as 
\begin{equation}
    (a,b) = \sum_{\omega \in W_A}a_\omega b_\omega
\end{equation}
where $a_\omega, b_\omega$ denote the coefficients in front of the word $\omega$ in $a,b$ respectively. Note that this sum is finite because $a$ is a finite linear combination of words. With this pairing, $\mathcal{A}^\infty$ can be identified as the algebraic dual space of $\mathcal{A}$. When restricted to $\mathcal{A} \times \mathcal{A}$, this pairing yields a scalar product with basis $W_A$ and dual basis $W_A'$. In the sequel we allow implicit and free conversion of letters and words, including the empty word $e$, according to context use the same notation $W_A$ for the word basis and its dual. 
\end{remark}






\begin{definition}
The \DEF{free magma} $\mathcal{M}_A$ is the minimal non-empty set satisfying: i)  $A \subset \mathcal{M}_A$, and ii) if $t', t'' \in \mathcal{M}_A$ then $(t', t'' )\in \mathcal{M}_A$. The \DEF{degree} of $t$ is defined recursively as $|t|=1$ if $t \in A$, otherwise if $t', t'' \in \mathcal{M}_A$ then $|t| = |t'| + |t''|$.
\end{definition}

\begin{remark}
\label{rem:abuse}
Let $V$ be a vector space.
The space $\mathcal B$ of bilinear maps $V\times V\to V$
naturally forms a magma, via composition.
For a fixed bilinear map $\phi: V\times V \to V$
and a set map $\iota: A \to V$ we abuse notation
and also write $\phi: \mathcal{M}_A \to \mathcal B$ for
the unique morphism of magmas characterized by
\begin{align*}
    \phi(a)        &= \iota(a), a \in A \\
    \phi((t',t'')) &= \phi( \phi(t'), \phi(t'') ).
\end{align*}
\end{remark}

\begin{definition}
The \DEF{foliage map} $f:\mathcal{M}_A \to W_A$ is defined on a letter $a \in A$ as $f(a)=a$ and on a tree $t=(t_1,t_2) \in \mathcal{M}_A$ as $f(t)=f(t_1)f(t_2)$ where the product is the tensor product (or concatenation of words). 
\end{definition}  

\begin{remark}
As noted in \citep{reutenauer1993london}, $\mathcal{M}_A$ can be equivalently identified with the set of binary, planar, rooted trees with leaves labelled in $A$. For a given element $t \in \mathcal{M}_A$ we will refer to the collection of letters appearing in its leaves as its foliage.
\end{remark}




\section{The free half shuffle algebra of Sch\"utzenberger}\label{sec:identities}

In this section we follow \cite{schutzenberger1958propriete} to define the half shuffle product and introduce the corresponding free algebra. We also provide two, to our knowledge, new identities in arity $3$ involving the commutator of the half shuffle product. These identities will be used in the next section to prove one of the main result of this paper.

\begin{definition}[\cite{schutzenberger1958propriete}]\label{def:half_shuffle}
The \DEF{(left) half shuffle product} $\lhs : \mathcal{A} \times \mathcal{A} \to \mathcal{A}$ is a bilinear form defined by extending uniquely, by linearity on the decomposition (\ref{eqn:z_decomposition}), the following relations
\begin{enumerate}
    \item $e \lhs f = 0 \lhs f = f \lhs 0 = 0 \ $ and $ \ f \lhs e = f, \quad $ for any $f \in \Anonempty, \quad $ 
    
    and by induction
    \item $f \prec f'' = a(f'\prec f'' + f'' \prec f'), \quad$ for any $f=af'$, with $a \in A, f' \in \Anonempty$, and $f'' \in \mathcal{A}$.
\end{enumerate}
\end{definition}

Note that the above definition of $\prec$ is independent of the choice of basis of $\mathcal{A}$.

\begin{remark}
    \Cref{def:half_shuffle} differs slightly from the usual algebraic convention that chooses to not define $e\lhs e$, as seen e.g. in \cite{ebrahimifardpatras2015}. In this paper, we follow to the letter \cite{schutzenberger1958propriete} where the half shuffle product is defined on $\Anonempty$ and $\langle e \rangle$, and then extended uniquely to a bilinear map on the direct sum (\ref{eqn:z_decomposition}) of these two spaces, that is to say the full algebra $\mathcal{A}$. Sch\"utzenberger refers to this canonical extension as \emph{prolongment}. 
\end{remark}

The following theorem is one the main results in   \cite{schutzenberger1958propriete}. 

\begin{theorem}\label{thm:schutzenberger}
    $\mathcal{A}$ is the free algebra over $A$ with respect to the half shuffle product $\prec$. 
\end{theorem}

We refer to this algebra as the \emph{free half shuffle algebra of Sch\"utzenberger}.


The \DEF{shuffle product} $\shuffle : \mathcal{A} \times \mathcal{A} \to \mathcal{A}$ is defined for any $f, g \in \mathcal{A}$ from the half shuffle $\prec$ as
\begin{equation}\label{eqn:shuffle}
    f \shuffle g = f \lhs g + g \lhs f + \langle f,e\rangle \langle g,e\rangle e.
\end{equation}
The algebra $(\mathcal{A},\shuffle)$ is an associative and commutative algebra known as the \DEF{shuffle algebra}. 
\begin{remark}
    Note that if $f,g \in \Anonempty$ then (\ref{eqn:shuffle}) reduces to the more conventional relation
    \begin{equation*}
        f \shuffle g = f \lhs g + g \lhs f.
    \end{equation*}
\end{remark}

The $\area$ operator is defined as the commutator of the half shuffle product and will be a core component of the main result in the next section.

\begin{definition}
    The operator $\area : \mathcal{A} \times \mathcal{A} \to \mathcal{A}$ is the bilinear form defined for $f, g \in \mathcal{A}$ as
    \begin{equation}
        \area(f, g) = f \lhs g - f \lhs g.
    \end{equation}
\end{definition}

In the next section we will provide concrete examples to demonstrate how Sch\"utzenberger's definition of half shuffle is completely consistent with classical integration on paths.

\subsection{Sch\"utzenberger's half shuffle is consistent with calculus}

Consider a smooth path $\gamma : [0,1] \to \mathbb{R}^d$, an interval $[a,b] \subset [0,1]$ and three elements $f,g,h \in \mathcal{A}$.  

Define the following one-dimensional paths on $[a,b]$:
\begin{equation*}
    \boldsymbol 1 : t \mapsto \langle e,\mathcal{S}(\gamma)_{a,t}\rangle, \quad f^\gamma : t \mapsto \langle f,\mathcal{S}(\gamma)_{a,t} \rangle, \quad  g^\gamma : t \mapsto \langle g,\mathcal{S}(\gamma)_{a,t}\rangle, \quad  h^\gamma : t \mapsto \langle h,\mathcal{S}(\gamma)_{a,t}\rangle.
\end{equation*}
Note that the path $\boldsymbol 1 \equiv 1$ is constantly equal to $1$. 

Notice how the relation $e \lhs f = 0$ in \Cref{def:half_shuffle} is consistent with the basic fact
\begin{align*}
    \langle e \prec f, \mathcal{S}(\gamma)_{a,t}\rangle = \int_a^t f^\gamma_s \D \boldsymbol 1_s = 0 = \langle 0, \mathcal{S}(\gamma)_{a,t}\rangle,
\end{align*}
 while the relation $f \lhs e = f - \langle f, e \rangle e$ is consistent with the fundamental theorem of calculus
\begin{align*}
  \langle f \prec e, \mathcal{S}(\gamma)_{a,t}\rangle = \int_a^t \boldsymbol 1_s \D f^\gamma_s = \int_a^t \D f^\gamma_s = f^\gamma_t - f^\gamma_a = \langle f - \left\langle f, e\rangle e, \mathcal{S}(\gamma)_{a,t}\right\rangle.
\end{align*}

All other classical rules of calculus follows. For example integration by parts 
\begin{equation*}
    \langle f \shuffle g - \langle f,e\rangle \langle g,e\rangle e, \mathcal{S}(\gamma)_{a,t} \rangle = f^\gamma_t g^\gamma_t - f^\gamma_a g^\gamma_a =  \int_a^t f^\gamma_s \D g^\gamma_s + \int_a^t g^\gamma_s \D f^\gamma_s = \langle f \lhs g + g \lhs f,\mathcal{S}(\gamma)_{a,t} \rangle ,
\end{equation*}
follows from the definition of shuffle product in equation (\ref{eqn:shuffle}).

Another classical example is provided by chain rule reads
\begin{equation*}
    \langle  f \lhs (g \shuffle h), \mathcal{S}(\gamma)_{a,t} \rangle = \int_a^t f^\gamma_s g^\gamma_s \D h^\gamma_s = \int_a^t f^\gamma_sd\left(\int_a^s g^\gamma_u \D h^\gamma_u\right) = \langle (f \lhs g) \lhs h, \mathcal{S}(\gamma)_{a,t} \rangle,
\end{equation*}
which matches the algebraic relation
\begin{equation}\label{eqn:chain-rule}
    f \lhs (g \shuffle h) = (f \lhs g) \lhs h. 
\end{equation}
Equation (\ref{eqn:chain-rule}) can be easily verified to hold for letters, and hence for all elements of $\mathcal{A}$ by freeness.



Next we present known and, to our knowledge, new identities on $\mathcal{A}$ involving $\lhs, \shuffle$ and $\area$.

\subsection{Identities}

The first identity is a direct application of the chain rule and integration by parts. When restricted to $\Anonempty$ it is known in the literature as \emph{Zinbiel identity} \cite{dzhumadil2007zinbiel}.

\begin{lemma}\label{lemma:zinbiel_identity}
For any $f,g,h \in \mathcal{A}$ the following identity holds
\begin{equation}\label{eq:mod_zinbiel}
   (f \lhs g) \lhs h = f \lhs (g \lhs h) + f \lhs (h \lhs g) + \langle g, e\rangle \langle h, e\rangle f \lhs e.
\end{equation}
\end{lemma}

\begin{proof}
A direct application of the chain rule and integration by parts yields
\begin{align*}
    (f \lhs g) \lhs h &= f \lhs (g \shuffle h)\\
    &= f \lhs (g \lhs h + h \lhs g  +  \langle g, e\rangle \langle h, e\rangle e)\\
    &= f \lhs (g \lhs h) + f \lhs (h \lhs g) +  \langle g, e\rangle \langle h, e\rangle (f - \langle f, e\rangle e),
\end{align*}
and the result follows from equation (\ref{eqn:shuffle}).
\end{proof}

\begin{remark}
When $f,g,h \in \Anonempty$ equation (\ref{eq:mod_zinbiel}) reduces to the Zinbiel identity
\begin{equation*}
   (f \lhs g) \lhs h = f \lhs (g \lhs h) + f \lhs (h \lhs g).
\end{equation*}
\end{remark}




\begin{remark}\label{lemma:shuffles_iter}
Using \Cref{lemma:zinbiel_identity} it is possible to obtain the following identity 
\begin{equation*}
    f_1\shuffle ... \shuffle f_n = \sum_{\sigma\in \mathfrak{S}_n}
     (... (f_{\sigma(1)}\lhs f_{\sigma(2)})\lhs ...)\lhs f_{\sigma(n)}
\end{equation*}
for any $n\geq 2$ and $f_1,...,f_n \in \Anonempty$, where $\mathfrak{S}_n$ is the symmetric group of order $n$.
\end{remark}




\begin{remark}
We note an important result obtained by \cite{dzhumadil2007zinbiel} stating that the $\area$ operator satisfies no further identity in arity three, but it does satisfy the so-called \emph{Tortkara identity} in arity four. While the Tortkara identity will play no further role in this paper, we mention it here for completeness: for any $f,g,h,i \in \Anonempty$, we equivalently have
\begin{equation*}
    \area(\area(f,g),\area(f,h))=\area(f,\vol(f,g,h))
\end{equation*}
and
\begin{align*}
    \area(\area(f,g),\area(i,h)) &+ \area(\area(h,g),\area(i,f))\\ &=\area(f,\vol(g,h,i))+\area(h,\vol(g,f,i))
\end{align*}
where $\vol(f,g,h):=\area(\area(f,g),h)+\area(\area(g,h),f)+\area(\area(h,f),g)$.

We furthermore note that Tortkara algebras have been studied more in \citep{dzhumadil2019speciality}, where it has been shown that the span inside $\mathcal{A}$ of iterated areas of letters forms a free Tortkara algebra for $|A|=2$, while the question remains open for larger alphabets.
\end{remark}

\begin{remark}[\textbf{left/right areas}]
In this paper, area is defined as the commutator of the left half shuffle. In \citep{diehl2020areas}, the right half shuffle is introduced and area is defined as the as the commutator of the right half shuffle. Although closely connected, these are not identical. The left half shuffle is consistent with \citep{reutenauer1993london} and matches the conventions for Hall basis used there (see later sections). The right half shuffle is more consistent with the convention used in integration as the integrand is on the left and the integrator is on the right. The reversed order of terms within equation (\ref{eqn:chain-rule}) reflects this dissonance. The proofs of our main results imply equivalent results with the other definition of area, by reversing everything. 
\end{remark}

Contrary to the Lie bracket $[\cdot,\cdot]$, $\area$ does not satisfy the Jacobi identity. However, it satisfies the following two non-trivial and, to our knowledge, new identities that will be leveraged to prove one of the main results of this paper in the next section. 

\begin{lemma}[\textbf{shuffle-pullout identity}]\label{lemma:shuffle-pullout}
    For any $f,g,h \in \mathcal{A}$ the following relation holds
    \begin{align*}
        3 \area(h, f \shuffle g) &= f \shuffle \area(h,g) + g \shuffle \area(h,f) - f \shuffle g \shuffle h + \langle f, e \rangle\langle g, e \rangle \langle h, e \rangle e\\  
        & + \area(\area(h,g),f) + \area(\area(h, f),g).
    \end{align*}
\end{lemma}

\begin{proof}
    It's easy to check that the relation holds for the empty word $e$ and for letters $a,b,c \in A$
    \begin{align*}
        3 \area(c,a\shuffle b) 
        &=-3\,abc-3\,acb-3\,bac-3\,bca+3\,cab+3\,cba\\
        &=a \shuffle \area(c,b) + b \shuffle \area(c,a) - a \shuffle b \shuffle c \\
        & + \area(\area(c,b),a) + \area(\area(c, a),b).
    \end{align*}
    By Theorem \ref{thm:schutzenberger} we know that $\mathcal{A}$ is free, as a half shuffle algebra over $A$, therefore the above relation extends to any triple of elements in $\mathcal{A}$.
\end{proof}

\begin{remark}
    When $f,g,h \in \Anonempty$ the shuffle-pullout identity in \Cref{lemma:shuffle-pullout} reduces to 
     \begin{align*}
        3 \area(h, f \shuffle g) &= f \shuffle \area(h,g) + g \shuffle \area(h,f) - f \shuffle g \shuffle h \\  
        & + \area(\area(h,g),f) + \area(\area(h, f),g) .
    \end{align*}
\end{remark}


\begin{lemma}[\textbf{area-Jacobi identity}]\label{lemma:area-jacobi}
For any triple $f, g, h \in \mathcal{A}$ the following relation is satisfied
    \begin{align*}
        &\area(\area(f, g),h) + \area(\area(g, h),f) + \area(\area(h, f), g) \\
         & =-f\shuffle \area(g, h) - g \shuffle \area(h,f) - h \shuffle \area(f, g).
    \end{align*}
\end{lemma}

\begin{proof}
     As before, the relation can be easily verified to hold for $e$ and for letters $a, b, c \in A$:
     \begin{align*}
        &\area(\area(a, b),c) + \area(\area(b, c),a) + \area(\area(c, a), b) \\
        &= -abc+acb+bac-bca-cab+cba\\
        & = -a \shuffle \area(b, c) - b \shuffle \area(c, a) - c \shuffle \area(a, b).
    \end{align*}
\end{proof}

\begin{remark}
 On $\Anonempty$, starting only from the identities 1) $f\shuffle g=g\shuffle f$, 2) $\area(f,g)=-\area(g,f)$, 3)  shuffle-pullout,
4) area-Jacobi, it follows from simple calculations that one can recover associativity for $\shuffle$ and the (left) Zinbiel identity for the left half shuffle $\lhs$, now defined by
$f \lhs g:=\frac{1}{2}(f\shuffle g+\area(f,g))$. Through the Zinbiel identity one then can show the Tortkara identity for $\area(f,g)=f\lhs g-g\lhs f$ as usual.



\end{remark}

\section{Polynomials in iterated areas}\label{sec:iterated_areas}

In this section we present our first main result, namely that polynomial in iterated areas generate the free half-shuffle algebra. We note that this result already appears in \citep{diehl2020areas}, however our proof is significantly shorter and based on induction.

\subsection{Polynomials in iterated areas are a generating set}

Recalling \Cref{rem:abuse}, we extend $\area$ to $\mathcal{M}_A$.
\begin{definition}
$f \in \mathcal{A}$ is an \DEF{iterated area} if there exists a tree $t \in \mathcal{M}_A$ so that $f=\areaTree(t)$. 

A shuffle monomial of shuffle-degree $n$ is the shuffle product of $n$ iterated areas 
\begin{equation}
    A_1 \shuffle ... \shuffle A_n.
\end{equation} 
The empty monomial $e$ has shuffle-degree $0$. A shuffle polynomial of shuffle-degree $n$ is a non-degenerate linear combination of such shuffle monomials. Its shuffle-degree is the maximal shuffle-degree of the monomials in the expression.
\end{definition}  

The sequence defined in the following lemma will play a role in what follows.  

\begin{lemma}\label{lemma:coeff}
The sequence of negative rationals $\beta_k = -(k - 1)/(k + 1)$ with $k\geq 1$ is monotone decreasing to $-1$ and satisfies the following recursion
\begin{equation}
    \beta_1 = 0, \quad \beta_k = \frac{\beta_{k-1} -1 }{\beta_{k-1}+3}.
\end{equation}
\end{lemma}

Exploiting the identities we introduced in the previous section we give a short and direct proof of the main result in \citep{diehl2020areas}.

\begin{theorem}\cite[Corollary 5.6]{diehl2020areas}\label{thm:iterated_areas}
Any element in $(\mathcal{A},\prec)$ can be written as a shuffle polynomial in iterated areas $\{\areaTree(t)\ |\ t \in \mathcal{M}_A\}$.
\end{theorem}

Before reproving the theorem we establish the following fundamental re-writing rule that allows one to rewrite the area of a shuffle polynomial in iterated areas with a single iterated area as a new shuffle polynomial in iterated areas, and provides an explicit expression for the monomial of highest shuffle-degree. The proof will crucially depend on both lemmas \ref{lemma:shuffle-pullout}, \ref{lemma:area-jacobi}.

\begin{theorem}\label{lemma:iterated-areas}
For any $n\geq 1$ and any $n+1$ iterated areas $A_1, ..., A_n, A$, the following relation holds 
\begin{equation}\label{eqn:iterated-areas}
    \area(A,A_1 \shuffle ... \shuffle A_n) = \beta_n A \shuffle A_1 \shuffle ... \shuffle A_n  + Q.
\end{equation}
where $\beta_n =-(n - 1)/(n + 1) $, and  $Q$ is a shuffle polynomial in iterated areas of shuffle-degree at most $n$.
\end{theorem}

\begin{remark}
 Note that it remains an open problem whether $\alpha=\beta_n$ is the only real number such that 
 \begin{equation*}
  \area(a, a_1 \shuffle ... \shuffle a_n)- \alpha a \shuffle a_1 \shuffle ... \shuffle a_n
 \end{equation*}
 can be expressed as a shuffle polynomial in iterated areas of shuffle-degree at most $n$ for any letters $a,a_1,\dots, a_n$. 
 This question arises due to the fact that iterated areas do not \emph{freely} generate the shuffle algebra.
 However, for the example $n=2$, $\beta_2=-1/3$ is indeed the only such coefficient because the area-Jacobi identity is the only relation between iterated areas on level $3$.
\end{remark}

\begin{proof}
We prove the statement (\ref{eqn:iterated-areas}) by induction on $n$. If $n=1$ then the statement is trivially true, with $\beta_1 = 0$ and $Q=\area(A_1, A)$. 

Suppose the statement (\ref{eqn:iterated-areas}) holds for any $n < k$. Consider $k$ iterated areas $A_1, ..., A_k$ and an additional iterated area $A$. We recall that the shuffle product $\shuffle$ is associative and commutative on $\mathcal{A}$. By the shuffle-pullout identity we have
\begin{align*}
    3 \area(A,A_1 \shuffle ... \shuffle A_k) &= A_1 \shuffle \area(A,A_2 \shuffle ... \shuffle A_k) \\
    &+ A_2 \shuffle ... \shuffle A_k \shuffle \area(A,A_1)\\
    &- A_1 \shuffle ... \shuffle A_k \shuffle A\\
    &+ \area(\area(A,A_1), A_2 \shuffle ... \shuffle A_k)\\
    &+ \area(\area(A,A_2 \shuffle ... \shuffle A_k), A_1).
\end{align*}

By induction ($n=k-1$) we have that
\begin{align*}
A_1 \shuffle \area(A,A_2 \shuffle ... \shuffle A_k) &= A_1 \shuffle(\beta_{k-1} A \shuffle A_2 \shuffle ... \shuffle A_k + Q_1')\\
&=\beta_{k-1} A \shuffle A_1 \shuffle ... \shuffle A_k + A_1 \shuffle Q_1'
\end{align*}
where $A_1 \shuffle Q_1'$ is a shuffle-polynomial of shuffle-degree $k$. By definition $\area(A,A_1)$ is a iterated area and so
$$Q_2' = A_2 \shuffle ... \shuffle A_k \shuffle \area(A,A_1)$$ 
is a shuffle monomial of shuffle-degree $k$. Similarly, the induction hypothesis implies that
$$Q_3' = \area(\area(A,A_1), A_2 \shuffle ... \shuffle A_k)$$

is a shuffle-polynomial of shuffle-degree $k$, where $\hat Q_3'$ is a shuffle-polynomial of shuffle-degree $k-1$. Hence, $Q' = A_k \shuffle Q_1' + Q_2' + Q_3'$ is a shuffle polynomial of shuffle-degree k and  
\begin{align}\label{eqn:utils}
    3 \area(A, A_1 \shuffle ... \shuffle A_k) &= Q' + (\beta_{k-1}-1) A_1 \shuffle ... \shuffle A_k \shuffle A\\
    &+ \area(\area(A, A_2 \shuffle ... \shuffle A_k), A_1). \nonumber
\end{align}

It remains to consider the last term $\area(\area(A, A_2 \shuffle ... \shuffle A_k), A_1)$. 

By the area-Jacobi identity and the anticommutativity of $\area$ we can rewrite this term as follows
\begin{align}
    \area(\area(A, A_2 \shuffle ... \shuffle A_k), A_1) &= \area(\area(A_1,A_2 \shuffle ... \shuffle A_k), A) \nonumber\\
    &- \area(\area(A_1, A), A_2 \shuffle ... \shuffle A_k) \nonumber \\
    &+ A \shuffle \area(A_1,A_2 \shuffle ... \shuffle A_k) \nonumber \\
    &- A_2 \shuffle ... \shuffle A_k \shuffle \area(A_1, A) \nonumber\\
    &- A_1 \shuffle \area(A,A_2 \shuffle ... \shuffle A_k). \nonumber
\end{align}

Again, $\area(A_1,A)$ is a iterated area, and by induction the term 
$$Q_1'' = - \area(\area(A_1, A), A_2 \shuffle ... \shuffle A_k)$$ 

is a polynomial in iterated areas of shuffle-degree at most $k$. The term 
$$Q_2'' = -A_2 \shuffle ... \shuffle A_k \shuffle \area(A_1, A)$$ 

is clearly a monomial in iterated areas of shuffle-degree $k$. By induction we have that
$$P_1 = \area(A_1,A_2 \shuffle ... \shuffle A_k) = \beta_{k-1} A_1 \shuffle ... \shuffle A_k + P_1'$$ 

where $P_1'$ is a polynomial in iterated areas of shuffle-degree $k-1$. Similarly
$$P_2 = \area(A,A_2 \shuffle ... \shuffle A_k) = \beta_{k-1} A \shuffle A_2 \shuffle ... \shuffle A_k + P_2'$$ 

where $P_2'$ is a polynomial in iterated areas of shuffle-degree $k-1$. Therefore
$$Q_3''= A \shuffle \area(A_1,A_2 \shuffle ... \shuffle A_k) = \beta_{k-1} A \shuffle A_1 \shuffle ... A_k + A \shuffle P_1'.$$

Similarly 
$$Q_4'' = - A_1 \shuffle \area(A, A_2 \shuffle ... \shuffle A_k) = - \beta_{k-1} A \shuffle A_1 \shuffle ... A_k - A_1 \shuffle P_2'.$$

Combining terms we get a cancellation and degree reduction so that
\begin{align*}
    Q_3'' + Q_4'' &= \beta_{k-1} A_1 \shuffle ... A_k \shuffle A + A \shuffle P_1' - \beta_{k-1} A_1 \shuffle ... A_k \shuffle A - A_1 \shuffle P_2' \\
    &=  A \shuffle P_1' - A_1 \shuffle P_2' 
\end{align*}
is a polynmomial in iterated areas of shuffle-degree $k$. Setting $Q'' = Q_1'' + Q_2'' + Q_3'' + Q_4''$ (which is a polynomial in iterated areas of shuffle degree $k$) and substituting in equation (\ref{eqn:utils}) we get 
\begin{align}
    3 \area(A,A_1 \shuffle ... \shuffle A_k) &= (\beta_{k-1}-1) A_1 \shuffle ... \shuffle A_k \shuffle A  \label{eqn:iterated-areas-final} \\
    &+ \area(\area(A_1,A_2 \shuffle ... \shuffle A_k), A) + Q' + Q'' \nonumber\\ 
    &= (\beta_{k-1}-1) A_1 \shuffle ... \shuffle A_k \shuffle A + \area(P_1, A) + Q' + Q''  \nonumber
\end{align}

$P_1$ being a polynomial in iterated areas of shuffle-degree $k-1$, we have by induction that $\area(P_1,A)$ is a polynomial in iterated areas of shuffle-degree $k$. Hence, by construction $Q = Q' + Q'' + \area(P_1, A)$ is a polynomial in iterated areas of shuffle-degree $k$. Therefore equation (\ref{eqn:iterated-areas-final}) becomes
\begin{align*}
    3 \area(A,A_1 \shuffle ... \shuffle A_k) &= (\beta_{k-1}-1) A_1 \shuffle ... \shuffle A_k \shuffle A \\
    &- \beta_{k-1}\area(A,A_1 \shuffle ... \shuffle A_k) + Q.
\end{align*}

Rearranging the terms we get the following final expression
\begin{align*}
    \area(A,A_1 \shuffle ... \shuffle A_k) &= \frac{\beta_{k-1}-1}{\beta_{k-1}+3} A_1 \shuffle ... \shuffle A_k \shuffle A + \frac{1}{\beta_{k-1} + 3}Q.
\end{align*}

Setting $\beta_k = \frac{\beta_{k-1}-1}{\beta_{k-1}+3} $ and noting that $\beta_1=0$ the result follows from Lemma \ref{lemma:coeff}.
\end{proof}

\begin{proof}[Proof of Theorem \ref{thm:iterated_areas}]
Since linear combinations of polynomials are polynomials, it suffices to prove that words in $W_A$ are polynomial in iterated areas. We prove by induction that every word $w \in W_A$ of length $|w|=n$ can be expressed as polynomial in iterated areas of shuffle-degree $n$. The result is trivial for $n=0$. Let $n\geq 1$. We assume that $w$ is a word of length $n>0$ and that any word of length $<n$ can be written as a polynomial in iterated areas of the appropriate degree. 

Since $|w|>0$, $w$ can be written as follows
\begin{equation}
w = av = a \lhs v
\end{equation}

where $v \in W_A$ is of word of length $|v|=n-1$ and $a \in A \subset \mathcal{A}$ is a letter. Moreover for any elements of $\mathcal{A}$
\begin{align}\label{eqn:simple}
    a \lhs v &= \frac{1}{2} (\area(a, v) + a \shuffle v - (a,e)(v,e)e) \\ \label{eqn:simple1}
    &= \frac{1}{2} (\area(a, v) + a \shuffle v)
\end{align}
since  $a$ is a letter.
The length of the word $v$ in (\ref{eqn:simple}) is equal to $n-1$, so by induction it can be written as a polynomial in iterated areas of shuffle-degree $n-1$. Hence, the term $a \shuffle v$ is a shuffle polynomial in iterated areas of shuffle-degree $n$. By Theorem \ref{lemma:iterated-areas} the term $\area(a,v)$ is also a polynomial in iterated areas of shuffle-degree $n$, and so $w$ a polynomial in iterated areas of shuffle-degree $n$. This concludes the induction and the proof.
\end{proof}

\section{A structure theorem for streamed information}\label{sec:structure_theorem}


To present our structure theorem we will need to introduce the free Lie algebra $\mathcal{L}_A$ over $A$. 

\subsection{The free Lie algebra}

$(\mathcal{A},[\cdot,\cdot])$ is also a Lie algebra with Lie bracket $[x,y]=x \otimes y - y \otimes x$ for $x,y \in \mathcal{A}$.

\begin{definition}
Denote by $(\mathcal{L}_A, [\cdot,\cdot])$ the Lie algebra generated by $A$ in $\mathcal{A}$, i.e. the intersection of all Lie algebras in $\mathcal{A}$ containing $A$.
\end{definition}

\begin{lemma}\cite[Theorem 0.5]{reutenauer1993london}
$(\mathcal{L}_A, [\cdot,\cdot])$ is the free Lie algebra over $A$.
\end{lemma}

\begin{remark}
The maps $\exp$ and $\log$ are classically defined as power series mapping $\mathcal{A}^\infty$ to $\mathcal{A}^\infty$. The truncated power series for $\exp^{(n)}$ and $\log^{(n)}$ provide good meaning for these operators as maps from $\mathcal{A}$ into $\mathcal{A}$. Those elements in $\mathcal{A}$ that are, at each truncated level $n\in\mathbb{N}$, in $\mathcal{L}_A$ are known as \DEF{Lie elements} and denoted by $\mathcal{L}^{(n)}_A$. Those elements in $\mathcal{A}$ that are, at each truncated level, exponentials of Lie elements, or equivalently, whose truncated logarithm is in $\mathcal{L}_A$, are known as \DEF{grouplike elements}  (and they form a group). The maps $\log$ and $\exp$ provide a one to one correspondence between group-like elements and Lie elements. 
\end{remark}

We report the following three classical results about the shuffle product $\shuffle$ and the free Lie algebra $\mathcal{L}_A$: the first states that the shuffle product characterises grouplike elements \cite[Lemma 2.17]{lyons2004differential}, the second provides a characterisation of Lie elements in $\mathcal{L}_A$ \cite[Theorem 3.1 (iv)]{reutenauer1993london}, and the third states that the exponential of Lie elements span the tensor algebra in a way that respects degrees of truncation
\cite[Lemma 3.4]{diehl2019invariants}.

\begin{theorem}\label{thm:three_statements_intro}
    Let $\ell \in \mathcal{L}_A$ be a Lie element. 
    \begin{enumerate}
        \item  $\langle f,\exp(\ell) \rangle \langle g,\exp(\ell)\rangle  =  \langle f \shuffle g,\exp(\ell)\rangle$ for any $f,g \in \mathcal{A}$.
        \item $\langle f \shuffle g, \ell\rangle = 0$ for any $f,g \in \Anonempty$.
        \item $\mathcal{A}^{(n)} = \text{Span}\{\exp^{(n)}(\ell) : \ell \in \mathcal{L}_A^{(n)}\}$ for any degree of truncation $n \in \mathbb{N}$.
    \end{enumerate}
\end{theorem}


In light of \Cref{thm:three_statements_intro} and of the following Lemma, the shuffle algebra $(\mathcal{A},\shuffle)$ can be identified with the algebra of $\mathbb{Q}$-polynomial functions on $\mathcal{L}_A$ with pointwise multiplication, denoted by $\mathbb{Q}[\mathcal{L}_A]$. 

\begin{lemma}\label{lemma:equivalence_poly}
For any $f \in \mathcal{A}$, the map $\ell \mapsto \langle f,\exp(\ell)\rangle $ is in $\mathbb{Q}[\mathcal{L}_A]$. Furthermore, the map $f \mapsto \langle f,\exp(\cdot)\rangle $ from $\mathcal{A}$ to $\mathbb{Q}[\mathcal{L}_A]$ is bijective. 
\end{lemma}

\begin{proof}
This result is classical, so we provide only a sketch of the proof. Any element $x \in \mathcal{A}$ is a finite sum of words in $W_A$ of some maximal length $d(x)$. Fix some basis $(\ell_i)_i$ for $\mathcal{L}_A$ that respects dimension and let $\ell = \sum l_i\ell_i$. Then the map $(s,\exp(\ell)) = (s,\exp(\sum_{d(\ell_i)\leq d(x)}l_i\ell_i))$ and the right hand side, truncated at degree $d(x)$ is clearly a polynomial in the $l_i$.  The exponentials of truncated Lie elements are linearly dense in the truncated tensor algebra, therefore $x$ is completely determined by its inner product with the $(x,\exp(\ell))$ as $\ell$ varies.
\end{proof}

\begin{remark}
It is an immediate corollary of these results, and of the \emph{Stone Weierstrass Theorem}, that any finite collection of distinct grouplike elements form the vertices of a simplex, and therefore that there is a linear functional that is one on any one of the elements and zero on the others.
\end{remark}

\begin{remark}
An analogy can be drawn with the Fourier transform seen as a change of basis for signals from time to frequency domain that turns point-wise multiplication into convolution. In our case, we can view $(\mathcal{A},\shuffle)$ as polynomial functions on $\mathcal{L}_A$ with pointwise multiplication, or as an algebra spanned by words, with the shuffle product, depending on our viewpoint.
\end{remark}

Next we introcude a special subsets of Hall trees in $\mathcal{M}_A$ classically used to construct bases for $\mathcal{L}_A$. Recall \Cref{rem:abuse} stating that any binary operator defined on words over $A$ automatically extends to an operator acting on trees from the magma $\mathcal{M}_A$.
In particular, this extends the Lie bracket, the half shuffle $\lhsf$, and the operation $\area$, to maps from $\mathcal{M}_A$ to $\mathcal A$.

\subsection{Hall sets}



\begin{definition}
A total order $<$ on a subset $M$ of $\mathcal{M}_A$ is an \DEF{ancestral order} if for any tree $t=(t',t'')$ of degree $\geq 2$ one has $t<t''$.
\end{definition}

This definition of ancestral order makes other constructions more transparent. It is obvious that ancestral orders exist on any magma and their restrictions to a subset are also ancestral. 


\begin{definition}
A subset $H$ of $\mathcal{M}_A$
together with an order $<$ on $H$ is a \DEF{Hall set} if the following conditions hold
\begin{enumerate}
\item $<$ is an ancestral order on $H$;
\item $A \subset H$;
\item for any tree $h=(h_1,h_2)\in \mathcal{M}_A$ of degree $\geq 2$, $h \in H$ if and only if: 
\begin{enumerate}
    \item $h_1, h_2 \in H$ and $h_1 < h_2$
    \item either $h_1 \in A$ or $h_2 \leq h_1''$ where $h_1=(h_1',h_1'')$.
\end{enumerate}
\end{enumerate}
\end{definition}

We note that, since $<$ is assumed to be ancestral, point 3.a implies $h < h_2$, which is a condition needed in the
general definition of Hall sets.
As pointed out in \cite[Proposition 4.1]{reutenauer1993london} and the surrounding discussion, Hall sets exist, any ancestral order on the full magma leads in a canonical way to to a unique Hall set, and that Hall sets are \DEF{closed}, i.e. each subtree of a Hall tree is again a Hall tree.

\begin{example}
The Hall set $H$ set used in the \texttt{esig} package \citep{esig} is defined as follows: elements are ordered so that they respect degree, and for any equal-length Hall trees $h=(h_1,h_2), h'=(h_1',h_2')$ their order is defined recursively as follows: $h<h'$ if either $h_1<h_1'$ or $h_1=h_1'$ and $h_2<h_2'$.
\end{example}

\begin{example}
Consider a total order on letters in $A$ and suppose that words in $W_A$ are ordered alphabetically. A \emph{Lyndon word} on $W_A$ is a non-empty word such that for any factorisation $\omega = uv$ with $u,v \in W_A$ non-empty one has $\omega < v$. Then, the set of Lyndon words ordered alphabetically is a Hall set \cite[Theorem 5.1]{reutenauer1993london}.
\end{example}

\begin{example}
Let $H_0 = A$ and order it totally. Define $H_{n+1}$ as the set of trees of the form
$$h = (...((h_1,h_2),h_3),...,h_k)$$
where $k \geq 2$ and $h_1,...,h_k \in H_n$ with
$$h_1<h_2 \geq h_3 \geq ... \geq h_k.$$ 
Now order $H_{n+1}$ totally. Finally let $H = \cup_{n \geq 0} H_n$ and extend the order in $H_n$ to $H$ by the condition
$$h_1 = H_m, h_2 \in H_n, m<n \implies h_1 > h_2.$$
Then $H$ is a Hall set \cite[Theorem 5.7]{reutenauer1993london}.
\end{example}


\begin{lemma}\cite[Corollary 4.14]{reutenauer1993london}\label{lemma:dim_hall}
Let $A$ be an alphabet of $q$ letters. The number of Hall trees of degree $n$ is equal to 
\begin{align}
    \mathcal{D}_H = \frac{1}{n} \sum_{d|n}\mu(d)q^{n/d}
\end{align}
where $\mu$ is the M\"obius function.
\end{lemma}

\subsection{The Poincar\'e-Birkhoff-Witt basis and its dual}

The Jacobi identities are linear relations between
degree-three Lie brackets arising from associativity of the underlying group operation. They make the derivation of a basis for the free Lie algebra $\mathcal{L}_A$ a deep and classic challenge. 

\begin{theorem}\cite[Theorem 4.9 (i)]{reutenauer1993london}\label{thm:Hall}
For any Hall set $H$, the collection of elements $\{[h] : h \in H\}$ form a linear basis for the free Lie algebra $\mathcal{L}_A$.
\end{theorem}

This basis admits a canonical extension to a basis of the tensor algebra $(\mathcal{A},\otimes)$. 

\begin{theorem}\cite[Theorem 4.9]{reutenauer1993london}\label{thm:PBW_TA}
    The decreasing products
    \begin{equation}
        [h_1]^{\otimes k_1}\otimes ... \otimes [h_n]^{\otimes k_n}, \quad h_i \in H, \quad  h_1>...>h_n
    \end{equation}
    is a basis of the tensor algebra  $(\mathcal{A},\otimes)$. This basis is called the \DEF{Poincar\'e-Birkhoff-Witt (PBW) basis}.
\end{theorem}

\begin{definition}
A word $\omega \in W_A$ is called a Hall word if $\omega$ is the image of a Hall tree $h \in H$ by the foliage map, i.e. $\omega=f(h)$. 
\end{definition}

\begin{remark}
    The foliage map is injective when restricted to a Hall set $H$ and there are efficient algorithms for recovering the Hall tree from a Hall word. 
\end{remark}

\begin{lemma}\cite[Corollary 4.7]{reutenauer1993london}
Every word $\omega \in W_A$ can be written uniquely as a decreasing product of Hall words 
\begin{equation}\label{eqn:word_decomposition}
    \omega = f(h_1)^{\otimes k_1} \otimes ... \otimes f(h_n)^{\otimes k_n}, \quad h_i \in H, \quad h_1>...>h_n.
\end{equation}
\end{lemma} 

\begin{remark}
If $\omega \in W_A$ is a word decomposed into its unique decreasing product of Hall words according to equation (\ref{eqn:word_decomposition}), then  $P_\omega$ is the corresponding PBW basis element as per Theorem \ref{thm:PBW_TA}
\begin{equation*}
    P_\omega = [h_1]^{\otimes k_1} \otimes ... \otimes [h_n]^{\otimes k_n}, \quad h_i \in H, \quad h_1>...>h_n.
\end{equation*}
\end{remark}

$\{P_\omega\}_{\omega \in W_A}$ is thus an enumeration of the PBW basis indexed by words. The next theorem provides exact formulae for the dual basis to the PBW basis.


\begin{theorem}\cite[Theorem 5.3]{reutenauer1993london}\label{thm:reut}
The dual basis $\{S_\omega\}_{\omega \in W_A}$ to the PBW basis $\{P_\omega\}_{\omega \in W_A}$ has the following properties:
\begin{enumerate}
    \item If $e$ is the empty word then $S_e = e.$
    \item If $\omega=f(h_1)^{\otimes k_1} \otimes ...\otimes f(h_n)^{\otimes k_n}$ is the unique factorization of the word $\omega $ in a decreasing product of Hall trees $h_1>...> h_n \in H$, then 
    \begin{equation}\label{eqn:Sw1}  
    S_\omega = \frac{1}{k_1! ... k_n!}S_{f(h_1)}^{\shuffle k_1} \shuffle ... \shuffle S_{f(h_n)}^{\shuffle k_n}. 
    \end{equation}
    \item If $h \in H$, then the word $f(h)=av$ for some letter $a \in A$ and word $v \in W_A$; moreover
    \begin{equation}
        S_{f(h)} = a \otimes S_v. \label{eqn:aSv}
    \end{equation}
    \end{enumerate}
\end{theorem}


Theorem \ref{thm:reut} is an important result due to Sch\"utzenberger and it is the structure theorem mentioned in the introduction. However, in the next section we provide our version of this theorem (which agrees with the version in \citep{sussmann1986product} but with a completely different proof) which consists of a more explicit recursive formula for the dual PBW basis elements $\{S_\omega\}_{\omega \in W_A}$ and identify them as Hall integrals. We note that this result is reported without proof also in \citep{kawski1999chronological,gehrig2008hopf}.

\subsection{Polynomials in Hall integrals are a free generating set}\label{sec:dualpbw}



\begin{definition}
An element $x$ of $\mathcal{A}$ is called a \DEF{Hall integral} if it is the image under the operator $\lhsf: \mathcal{M}_A \to \mathcal{A}$ of a Hall tree. That is to say, there exists a Hall tree $h \in H \subset \mathcal{M}_A$ so that $x = \lhsf(h)$. A (shuffle) polynomial in Hall integrals is a sum of shuffle monomials in Hall integrals.
\end{definition}

The following Lemma follows immediately from the definition of a Hall tree.

\begin{lemma}
Any Hall tree $h \in H$ can be uniquely decomposed as 
\begin{equation}
    h = (h_1h_2^k) = (...((h_1,h_2),h_2),... h_2)
\end{equation}
with $h_1, h_2\in H$,
and either $h_1$ is a letter
or $h_1'' \not= h_2$
and where the $h_2$ bracketing is repeated $k$ times. This is often referred to as the \DEF{Lazard decomposition} of $h$.
\end{lemma}

\begin{definition}
If $h=(h_1h_2^k)$ is the Lazard decomposition of a Hall tree $h \in H$ then we define the \DEF{Lazard depth} $\alpha_h$ of $h$ to be $1/k$. The \DEF{accumulated Lazard depth} of a Hall tree $h \in H$ is defined recursively: $\mathcal{A}_h = 1$ if $h \in A$, otherwise $h = (h',h'')$ and $\mathcal{A}_h = \alpha_h \mathcal{A}_{h'} \mathcal{A}_{h''}$.
\end{definition}

The following are the main results of this section.

\begin{theorem}\label{thm:hs}
For any Hall tree $h \in H \setminus A$ one has $h=(h',h'')$ and 
\begin{equation}
    S_{f(h)} = \alpha_h  \left( S_{f(h')} \lhs S_{f(h'')} \right)
\end{equation}
where $\alpha_h \in \mathbb{Q}$ is the Lazard depth of $h$.
\end{theorem}

\begin{theorem}\label{thm:hs2}
    For any Hall tree $h \in H$ one has
\begin{equation}
    S_{f(h)} = \mathcal{A}_h ( \lhsf(h))
\end{equation}
where $\mathcal{A}_h \in \mathbb{Q}$ is the accumulated Lazard depth of $h$.
\end{theorem}

\begin{theorem}\label{thm:hs1}
Consider all decreasing sequences $ h_i \in H$,  $h_1>...>h_n$, and strictly positive integers $k_i > 0$; then the elements
\begin{equation}
    S_\omega = \frac{\mathcal{A}_{h_1}^{k_1} \dots \mathcal{A}_{h_n}^{k_n}}{k_1! ... k_n!}(\lhsf(h_1))^{\shuffle k_1} \shuffle ... \shuffle (\lhsf(h_n))^{\shuffle k_n} 
\end{equation}
are the dual basis in  $\mathcal{A}$ to the PBW basis  $\{P_\omega=[h_1]^{\otimes k_1} \otimes ... \otimes [h_n]^{\otimes k_n}\}_{\omega \in W_A}$. 

\end{theorem}

Before proving Theorem \ref{thm:hs} we need the following combinatorial lemma.

\begin{lemma}\cite[Corollary 5.14]{reutenauer1993london}\label{lemma:dec}
Let $h=(h',h'') \in H$ be a Hall tree. Now $f(h)=av$, where $a \in A$ and $v \in W_A$.  Let $v=f(h_1)^{\otimes k_1} \otimes ... \otimes f(h_n)^{\otimes k_n}$ be the unique factorization of the word $v$ in a decreasing product of Hall trees $h_1>...> h_n \in H$. Then 
\begin{equation}
    h''=h_n.
\end{equation}
\end{lemma}

\begin{proof}[Proof of Theorem \ref{thm:hs}]
We write $f(h)=av$, with $a\in A$ and $v \in W_A$. Let $v=f(h_1)^{\otimes k_1} \otimes ... \otimes f(h_n)^{\otimes k_n}$ be the unique factorization of the word $v$ in a decreasing product of Hall trees $h_1>...> h_n \in H$. By Lemma \ref{lemma:dec} $h''=h_n$. By Theorem \ref{thm:reut} we also know that
\begin{align}
S_{f(h)}&=a \otimes S_v \label{eqn:a}\\
&=S_a\lhs S_v\label{eqn:b}\\
&=\frac{1}{k_1! ... k_n!} S_a\lhs (S_{f(h_1)}^{\shuffle k_1}\shuffle ... \shuffle S_{f(h_n)}^{\shuffle k_n})\label{eqn:c}\\
&=\frac{1}{k_1! ... k_n!} S_a\lhs ((S_{f(h_1)}^{\shuffle k_1}\shuffle ... \shuffle S_{f(h_n)}^{\shuffle k_n-1}) \shuffle S_{f(h_n)})\label{eqn:d}\\
&=\frac{1}{k_1! ... k_n!} S_a\lhs ((S_{f(h_1)}^{\shuffle k_1}\shuffle ... \shuffle S_{f(h_n)}^{\shuffle k_n-1}) \shuffle S_{f(h'')})\label{eqn:e}\\
&=\frac{1}{k_1! ... k_n!} (S_a\lhs (S_{f(h_1)}^{\shuffle k_1}\shuffle ... \shuffle S_{f(h_n)}^{\shuffle k_n-1}))\lhs S_{f(h'').}\label{eqn:f}
\end{align}
Equation (\ref{eqn:a}) is a restatement of (\ref{eqn:aSv}) in Theorem \ref{thm:reut}. Equation (\ref{eqn:b}) is immediate from the definition of $\lhs$. Equation (\ref{eqn:c}) follows from (\ref{eqn:Sw}) in Theorem \ref{thm:reut}.  Equation (\ref{eqn:d}) is simply the associative property of shuffle. Equation (\ref{eqn:e}) follows from Lemma \ref{lemma:dec}.  Equation (\ref{eqn:f}) follows from the chain rule (\ref{eqn:chain-rule}). Note that the inner term in equation (\ref{eqn:f}) can be reinterpreted as $S_{f(h')}$ (up to scalar) because
\begin{align}
S_{f(h')}
&=\frac{1}{k_1! ... (k_n-1)!} S_a\lhs (S_{f(h_1)}^{\shuffle k_1}\shuffle ... \shuffle S_{f(h_n)}^{\shuffle k_n-1}).\label{eqn:i}
\end{align}
Substituting this into equation (\ref{eqn:f}) and recalling the definition of the Lazard depth $\alpha_h$ we obtain
\begin{align}
S_{f(h)}&=\frac{1}{k_n} (S_{f(h')}\lhs S_{f(h'')})\\
&=\alpha_h (S_{f(h')}\lhs S_{f(h'')}).
\end{align}
\end{proof}

\begin{proof}[Proof of Theorem \ref{thm:hs2}]
We may proceed by induction. For any Hall tree $h \in A$ one has $S_{f(h)} = h \in \mathcal{A}$, $\lhsf(h)=h$, and $\mathcal{A}_h=1$ and so the theorem is true. On the other hand if $h=(h',h'')$ then, assuming the result holds for $h'$, $h''$: 
\begin{align}
    S_{f(h)} &= \alpha_h  \left( S_{f(h')} \lhs S_{f(h'')} \right)\\
   &= \alpha_h  \left((\mathcal{A}_{h'} (\lhsf(h'))) \lhs (\mathcal{A}_{h''}(\lhsf(h''))) \right)\\
   &= \mathcal{A}_h (\lhsf(h))
\end{align}
where we use Theorem \ref{thm:hs} for the first step, the truth of the result for $h'$ and $h''$ for the second, and the recursive definitions of $\lhsf(h)$ and $\mathcal{A}_h$ for the third. So the result is true for $h$. 
\end{proof}

\begin{proof}[Proof of Theorem \ref{thm:hs1}]
Recall from Sch\"utzenberger's theorem (Theorem \ref{thm:reut} in this paper) that any element $S_w$ in the dual basis to the PWB basis can be expressed uniquely as a shuffle monomial in $S_{f(h)}$. More precisely, consider the unique factorization of the word $w$ as a decreasing product of Hall words $w=f(h_1)^{\otimes k_1} \otimes ... \otimes f(h_n)^{\otimes k_n}$  where $h_1>...> h_n \in H$, then the dual basis element
    \begin{equation}
    S_w = \frac{1}{k_1! ... k_n!}S_{f(h_1)}^{\shuffle k_1} \shuffle ... \shuffle S_{f(h_n)}^{\shuffle k_n} \in \mathcal{A}.\label{eqn:Sw}    
    \end{equation}
Theorem \ref{thm:hs2} allows for $i = 1 ... n$ the substitution of $\mathcal{A}_{h_i} (\lhsf({h_i}))$ for $S_{f({h_i})}$ in this formulae which gives the specified expression for the dual basis element in terms of Hall integrals.
\end{proof}

In this section we have provided formulae for the dual PBW basis elements alternative but equivalent to the ones to be found in the book \citep{reutenauer1993london}.

\subsection{A conjecture}

\Cref{thm:hs1} states that polynomials in \emph{Hall integrals} freely generate the half shuffle algebra $(\mathcal{A},\prec)$ as an associative and commutative algebra. A natural question is whether a similar structure theorem holds in the case where the half shuffle $\prec$ on Hall trees is replaced by the commutator $\area$ as basic operation. This question has been, and still remain, a conjecture well supported by calculation for the last decade. 

\paragraph{Conjecture} \emph{Any element of $\mathcal{A}$ can be written uniquely as a polynomial over Hall areas} $\{\areaTree(h)\}_{h \in H}.$

Trying to solve this conjecture led us to consider an argument related to the well-known Lazard's elimination \citep{reutenauer1993london} to construct a canonical, but to our knowledge, new decomposition of the half shuffle algebra $(\mathcal{A},\prec)$ as shuffle power series in the greatest letter $c$ of the alphabet $A$ with coefficients in a sub-algebra $\mathcal{X}$ freely generated by a new alphabet $X$ with an infinite number of letters defined in terms of $c$ and all other letters in $A$. This construction, that we refer to as elimination trick, allows us to provide, in the next section, a second proof relying on an induction argument of our structure theorem.

\subsection{Another proof of the structure theorem}

The following simple and concrete observation will be expanded in this section.

If $(\mathcal{A},\prec)$ is the free half shuffle algebra over $A$, and $c\in A$ , and $X$ is the subset of $\mathcal{A}$ comprising $\frac{1}{k}\lhsf((ac^k))$, $a\in A\setminus c$, and $Z$ is the space spanned by words that do not begin with $c$; then $(Z,\prec)$ is a half shuffle algebra generated by $X$ in $\mathcal{A}$; moreover, $Z$ is freely generated as a half shuffle algebra by $X$, and therefore canonically isomorphic as a half shuffle algebra to the free half shuffle algebra $(\mathcal{X},\prec)$ over $X$. In characteristic zero, $Z$ is the half shuffle sub-algebra of $\mathcal{A}$ spanned by the words that do not begin with $c$. It is complimentary to $\mathcal{A}\shuffle c$ and we have 
\begin{align}
\mathcal{A} &= Z \oplus (\mathcal{A}\shuffle c) \nonumber\\
     & = Z \oplus (Z\shuffle c)  \oplus (\mathcal{A}\shuffle c \shuffle c) \label{eqdirect_sum}\\
     &=\dots \nonumber
    \end{align}
and any element in $\mathcal{A}$ can be expressed canonically as a shuffle power series in $c$ with coefficients in the half shuffle subalgebra $Z$. One can repeat this process by choosing a letter $d\in X$, and expanding every coefficient as a power series in $d$ with coefficients in the half shuffle subalgebra generated by the elements $\{\frac{1}{k}\lhsf((ad^k)), a\in X\setminus d\}$. In what follows we will make this precise.

\begin{definition}
Let $c$ be the greatest element of $A$ with respect to an ancestral ordering $<$. Define the subset of trees
\begin{equation}
    X = \{(ac^n), a \in A \setminus\{c\}, n \geq 0\} \subset \mathcal{M}_A.
\end{equation}
\end{definition}

With this choice of (infinite) alphabet, the following spaces and operators are automatically defined in the same way as their $A$ counterparts:

\begin{itemize} 
    \item $\mathcal{M}_X$ the free magma;
    \item $W_X$ the space of words in the alphabet $X$;
    \item $\mathcal{X}$ the vector space spanned by words in $W_X$;
    \item $\otimes_X$, $[\cdot,\cdot]_X$, $\lhs_X$,  $\areaTree_X$, $(\cdot,\cdot)_X$ the products and pairing on these spaces;
    \item $\mathcal{L}_X$ the free Lie sub-algebra of $(\mathcal{X},\otimes_X)$;
    \item $\exp_X$,  $\log_X$ the tensor series for the respective maps.
\end{itemize}

\begin{remark}
Note that the elements of $W_X$ are words whose letters are particular words in $A$. 
\end{remark}

\begin{theorem}\cite[Theorem 0.6]{reutenauer1993london}\label{thm:elimination} The Lie algebra $\mathcal{L}_A$ is the semi-direct product of $\mathcal{L}_X$ and $\mathbb{R}c $ 
\begin{equation}
    \mathcal{L}_A =  \mathcal{L}_X \ltimes \mathbb{R}c.
\end{equation}
\end{theorem}

As a result of \Cref{thm:elimination}, $\mathcal{L}_X$ is a Lie ideal and sub-algebra of co-dimension one in $\mathcal{L}_A$ and in particular $\mathcal{L}_A =  \mathcal{L}_X \oplus \mathbb{R}c$.

Next we report an important lemma from \cite{reutenauer1993london} which provides a very simple relation between Hall sets in $\mathcal{M}_A$ and Hall sets in $\mathcal{M}_X$.

\begin{lemma}\cite[Lemma 4.19 \& Section 5.6.3]{reutenauer1993london}\label{lemma:Hall}
    The unique homomorphism of magmas $\phi:\mathcal{M}_X \to \mathcal{M}_A$ that sends $x=(ac^n) \in X$ to $(ac^n) \in \mathcal{M}_A$ is an injection of magmas and its range is the free magma over $X$.  Furthermore $\phi(H_X) = H \cap \phi(\mathcal{M}_X)$ is the Hall set in $\mathcal{M}_X$ associated with the ordering $<$, $H = \{c\} \cup \phi(H_X)$, and $c$ is the greatest element of $H$.
\end{lemma}

\begin{remark}
$\mathcal{M}_X$ is a sub-magma and  inherits an ancestral ordering from $\mathcal{M}_A$. It follows that the image by $\phi$ of the Hall set $H_X$ in $\mathcal{M}_X$ associated to the ordering $<$ is $H\cap \phi(\mathcal{M}_X)$ (Lemma \ref{lemma:Hall}).
\end{remark}

When switching back and forth between the $X$- and $A$-spaces, the first objects one needs to have control over are letters from the two alphabets $X$ and $A$. In the next lemma we express the images under the various operators discussed so far of letters in $X$, seen as trees in $\mathcal{M}_A$, in terms of words from $W_A$.

\begin{lemma}\label{lemma:simple}
For any $x\in X$, the image $\phi(x)$ in $\mathcal{M}_A$ is of the form $(ac^n)$ for some $a \in A$ and $n\geq0$. The image of $(ac^n)$ under the operators $[\ ]$, $\lhs$, $\areaTree$ in $\mathcal{A}$, expressed in terms of words in $W_A$ are given by
\begin{align}
    [\phi(x)] = [(ac^n)] &= \binom{n}{0} ac ... c - \binom{n}{1} c a c ... c + ... + (-1)^n \binom{n}{n} c ... ca \label{eqn:lie_ac}\\
    \lhsf(\phi(x)) &= \lhsf((ac^n)) = n! a c... c \label{eqn:hs_ac}\\
    \areaTree(\phi(x)) &= \areaTree((ac^n)) = n!(a c... c - c a c ... c)  \label{eqn:area_ac}
\end{align}
where all the words are of length $n+1$ and contain exactly once the letter $a$.
\end{lemma}

The proof is left as an exercise to the reader.

The next lemma tells the relationship between integrals and areas on letters from $X$.

\begin{lemma}\label{lemma:simple_relation}
For any tree $(ac^n) \in \mathcal{M}_A$ one has
\begin{equation}\label{eqn:prec_acn}
    \lhsf((ac^n)) = \frac{1}{n+1} \areaTree((ac^n)) + \frac{n}{n+1} c \ \shuffle \lhsf((ac^{n-1})).
\end{equation}
\end{lemma}

\begin{proof}
From Lemma \ref{lemma:simple} we deduce the following identity
\begin{equation*}
    \areaTree((ac^n)) + (c \shuffle n \lhsf((ac^{n-1})))= (n+1) \lhsf((ac^n))
\end{equation*}
which after rearranging yields equation (\ref{eqn:prec_acn}).
\end{proof}

\begin{lemma}\label{lemma:simple_relation2}
For any $(ac^n) \in \mathcal{M}_A$ one has
\begin{equation}
    \lhsf((ac^n)) = \frac{1}{n+1} \sum_{k=0}^n c^{\shuffle k} \shuffle \areaTree((ac^{n-k})).
\end{equation}
\end{lemma}

\begin{proof}
This follows immediately from Lemma \ref{lemma:simple_relation} and an induction on $n$.
\end{proof}

\begin{remark}
Recall that the Lie bracket operator $[\cdot]$ is defined on $\mathcal{M}_A$ with values in $\mathcal{L}_A$. The restriction of $[\cdot]$ defined on $\mathcal{M}_A$ to $\mathcal{M}_X$ agrees with the natural definition of $[\cdot]$ on $\mathcal{M}_X$. It is also a simple exercise to prove that this compatibility between the restriction and the intrinsically defined operators holds for the tensor product and the Lie bracket.
\end{remark}

\begin{definition}
We denote by $J_c : \mathcal{X} \to \mathcal{A}$ the  unique $\lhs$-homomorphism that, by freeness of $(\mathcal{X},\prec_X)$ over $X$, extends to $\mathcal{X}$ the map 
\begin{equation}\label{eqn:J_c_x}
    (ac^n) \mapsto \frac{1}{n} (\lhsf((ac^n))), \quad n>0.
\end{equation}
Denote by $(Z,\prec)$ the half shuffle subalgebra of $(\mathcal{A},\prec)$ generated by the elements 
$$\{J_c(x) : x \in X\}.$$
\end{definition}

Next we prove that that the algebra $(Z,\prec)$ is closed under $\lhs$ and provide a characterisation of $Z$ as the linear span of words in $W_A$ that do not begin with the letter $c$.

\begin{lemma}\label{lemma:Z}
$Z$ is the span of words in $W_A$ that do not begin with the letter $c$
\begin{equation*}
    Z = \text{Span}\{w=a \prec v \in W_A \mid a \neq c, a \in A, v \in W_A\}.
\end{equation*}
In particular $Z$ is closed under $\lhs$.
\end{lemma}

\begin{proof}
Let $Z'$ be the linear span in $\mathcal{A}$ of the $w\neq cv \in W_A$. It is immediate from the definitions of $\shuffle$ and $\lhs$ on words that $Z'$ is closed under both operations. If $t=(t_1,t_2) \in \mathcal{M}_A$ and if, for $i=1,2$, 
$\lhsf(\phi(t_i)))\in Z'$ then  $\lhsf(\phi(t)) = (\lhsf(\phi(t_1))) \lhs (\lhsf(\phi(t_2)))$ is also in $Z'$ because $Z'$ is closed under $\lhs$. Let $x \in X$, then by equation (\ref{eqn:hs_ac}) 
\begin{align}
    \lhsf(\phi(x)) &= \lhsf((ac^n)) = n! a c... c \in Z'.
\end{align}
We may proceed recursively to see that every $\lhsf(t)$ contained in $Z$ is also an element of $Z'$; since $Z$ is generated by $\{J_c(x) : x \in X\}$ we conclude that $Z\subset Z'$. The unique decomposition of words into decreasing sequences of Hall words shows that the dimension of $Z'$ and $Z$ are equal, hence $Z=Z'$. 
\end{proof}

\begin{lemma}\label{lemma:decomposition}
The half shuffle algebra $\mathcal{A}$ has the following decomposition
\begin{equation}\label{eqn:decomposition}
    \mathcal{A} = Z \oplus (Z\shuffle c) \oplus  (Z\shuffle c^{\shuffle 2}) + ...
\end{equation}
\end{lemma}

\begin{proof}
Consider any word $w \in W_A$ beginning with $n$ number of $c$'s.
\begin{equation*}
    w = c \prec (c \prec( ...\prec (c \prec v)...))
\end{equation*} 
where $v=av' \in Z$ is a word that doesn't begin with $c$, i.e. $a \in A$, $a \neq c$, $v' \in W_A$. If $n=0$ then $w \in Z$. By induction on $n$
\begin{equation*}
    c \prec (c \prec( ...\prec (c \prec v)...)) - \alpha_n c^{\shuffle n} \shuffle v = L
\end{equation*}
where $\alpha_n \in \mathbb{R}$ and $L$ is a linear combination of words that begin with $k<n$ number of $c$'s. Hence, by induction on the number of $c$'s in front of the words, the word $w$ can be written as a shuffle polynomial in $c$ with coefficients in $Z$.
\end{proof}

\begin{lemma}\label{lemma:image_hall_integrals}
$J_c$ maps polynomials in Hall integrals $\lhs_X(h)$, $h \in H_X$ to polynomials in Hall intgrals $\lhs(\phi(h))$.
\end{lemma}

\begin{proof}
This follows immediately because $J_c$ is a half shuffle (and so shuffle) homomorphism.
\end{proof}

We now repeat our structure theorem and provide an alternative proof based on the elimination trick discussed so far in this section.

\begin{theorem}\label{thm:dual_pbw}
 The half shuffle algebra $(\mathcal{A},\prec)$ is freely generated by polynomials in Hall integrals $\lhs(h)$ for $h \in H$.
\end{theorem}

\begin{proof}
We can assume by induction that the theorem holds for $\mathcal{X}$, i.e. that $(\mathcal{X},\prec_X)$ is freely generated by polynomials in $\lhs_X(h)$ for $h \in H_X$. By Lemma \ref{lemma:image_hall_integrals}, $Z$ is freely generate by polynomials in $\lhs(h)$ with $h \in \phi(H_X)$. By \Cref{lemma:Hall}, $H = \{c\}\cup\phi(H_X)$ and with the decomposition (\ref{eqn:decomposition}) we conclude that $\mathcal{A}$ is freely generated by polynomials in $\lhs(h)$, $h \in H$.
\end{proof}

\subsection{Scalable computations of path signatures}

As mentioned in the introduction, instances of streamed information can be represented as a path $\gamma : [0,1] \to V$ with values on some finite dimensional vector space $V \simeq \mathbb{R}^d$, such a path is faithfully represented, up to reparameterisation, by the signature $\mathcal{S}_\gamma \in (\mathcal{A},\otimes)$.
Furthermore,
since the extended tensor algebra $(\mathcal{A},\otimes)$
is the algebraic dual of the half shuffle algebra $(\mathcal{A},\prec)$,
it is automatic to see that the restriction of linear functionals on $\mathcal{A}$ to the range of the signature form a unital algebra of real-valued functions that separates signatures. Hence, by the Stone-Weierstrass theorem linear functionals acting on the signatures are dense in the space of continuous, real-valued functions on compact sets of unparameterised paths. Thus, non-linear regression on pathspace can be realised by linear regression on the terms of the signature. However, terms in the signature contain some redundancy, which represents a major scalability issue, particularly because the number of distinct and linearly independent iterated integrals grows exponentially in the truncation level. In this paper, and in particular in Theorems \ref{thm:hs1} and \ref{thm:dual_pbw},
we identified sets of Hall integrals that can be used to compute any term in the signature with a minimal amount of computations. 

To illustrate this we consider a simple example. Let $d=3$ and let us identify the $3$-dimensional vector space $V$ as the space spanned by an alphabet of three letters $A = \{1,2,3\}$. Let $\omega = \boldsymbol{233212222111}$; note that $|\omega| = 12$. Then, computing the coefficient $(S_\omega,\mathcal{S}_\gamma)$ in the signature using existing software \citep{signatory, esig, reizenstein2018iisignature} (based on the \DEF{Chen's relation}) involve evaluating the level-$12$ truncated tensor exponential of increments $\exp^{(12)}(\gamma_t - \gamma_s)$. This operation has space and time complexities of $\mathcal{O}(3^{12})$. 

Instead, considering for example the Lyndon basis, one can precompute the factorisation of $\omega$ into decreasing product of Lyndon words and find
$$\omega = f(\boldsymbol 1)^{\otimes 3} \otimes f(((((\boldsymbol 1, \boldsymbol 2), \boldsymbol 2), \boldsymbol 2), \boldsymbol 2)) \otimes f(\boldsymbol 2) \otimes f(((\boldsymbol 2,\boldsymbol 3), \boldsymbol 3))$$
Therefore, by \Cref{thm:hs1} one has
$$S_\omega = \lhsf(\boldsymbol 1)^{\shuffle 3} \shuffle \frac{1}{4!} \lhsf(((((\boldsymbol 1, \boldsymbol 2), \boldsymbol 2), \boldsymbol 2), \boldsymbol 2)) \shuffle \lhsf(\boldsymbol 2) \shuffle \frac{1}{2!}\lhsf(((\boldsymbol 2,\boldsymbol 3), \boldsymbol 3)).$$
Using the interplay between algebraic operations $\prec$ and $\shuffle$ and the rules of calculus on paths outlined in \Cref{sec:identities} we obtain
\begin{align*}
    (S_\omega,\mathcal{S}_\gamma) &= \frac{1}{48} \alpha_1\alpha_2\alpha_3\alpha_4
\end{align*}
where
\begin{align*}
    \alpha_1 &= \int_0^1 d \gamma^{(1)}_t\\
    \alpha_2 &= \int_0^1\left(\int_0^v \left(\int_0^u\left(\int_0^t\gamma^{(1)}_sd\gamma^{(2)}_s\right)d\gamma^{(2)}_t\right) d\gamma^{(2)}_u\right)d\gamma^{(2)}_v\\
    \alpha_3 &= \int_0^1 d \gamma^{(2)}_t\\
    \alpha_4 &= \int_0^1\left(\int_0^{t} \gamma^{(2)}_{s} d\gamma^{(3)}_{s}\right)d\gamma^{(3)}_{t}.\\
\end{align*}

\section{Conclusion}

In this paper, we identified the free Zinbiel algebra introduced by \cite{schutzenberger1958propriete} with an algebra of real-valued functions on paths. We provided two, to our knowledge, new basic identities in arity 3 involving its symmetrization $\shuffle$ and its anti-symmetrization $\area$.
We showed that these are sufficient to recover the Zinbiel and Tortkara identities introduced by \citet{dzhumadil2007zinbiel}. We then used these identities to provide a direct proof of the main result in \citep{diehl2020areas} stating that polynomials in iterated areas generate the free Zinbiel algebra \citep{sussmann1986product}. Subsequently, we introduced minimal sets of Hall integrals and showed, with two different proof techniques, that polynomial functions on these Hall integrals freely generate the half shuffle algebra. This result can be interpreted as a structure theorem for streamed information, allowing to split real valued functions on streamed data into two parts: a first that extracts and packages the streamed information into Hall integrals, and a second that evaluates a polynomial in these without further reference to the original stream.
    

\section*{Acknowledgments}
We deeply thank Prof. Pavel Kolesnikov and Prof. Fr\'ed\'eric Patras for the helpful discussions and suggestions. Terry Lyons and C. Salvi’s contributions to this work is supported by the EPSRC program grant DataSig [grant number EP/S026347/1]. Terry Lyons' contribution was also supported by project partners: in part by The Alan Turing Institute under the EPSRC grant EP/N510129/1, in part by The Alan Turing Institute’s Data Centric Engineering Programme under the Lloyd’s Register Foundation grant G0095, in part by The Alan Turing Institute’s Defence and Security Programme, funded by the UK Government, and in part by the Hong Kong Innovation and Technology Commission (InnoHK Project CIMDA).

\bibliographystyle{apalike} 
\bibliography{references}

\end{document}